\documentclass[11pt,leqn]{amsart}
\usepackage{amssymb}
\usepackage{mathrsfs}
\usepackage[all]{xy}
\usepackage{color}
\usepackage{soul}
\usepackage{hyperref}
\usepackage{dsfont}

\newcommand\N{{\mathbb N}}
\newcommand\R{{\mathbb R}}
\newcommand\T{{\mathbb T}}
\newcommand\C{{\mathbb C}}

\newcommand\Z{{\mathbb Z}}

\def\CC{{\mathcal C}}

\def\HH{{\mathcal H}}

\def\LL{{\mathcal L}}

\def\eps{{\varepsilon}}

\newtheorem{theo}{Theorem}

\newtheorem{lem}[theo]{Lemma}

\newtheorem{rem}[theo]{Remark}

\newcommand{\beqn}{\begin{equation}}
\newcommand{\eeqn}{\end{equation}}
\newcommand{\bal}{\begin{aligned}}
\newcommand{\eal}{\end{aligned}}
\newcommand{\bear}{\begin{eqnarray}}
\newcommand{\eear}{\end{eqnarray}}
\newcommand{\bean}{\begin{eqnarray*}}
\newcommand{\eean}{\end{eqnarray*}}

\newcommand{\la}{\langle}
\newcommand{\ra}{\rangle}

\title[Landau damping for the linearized Vlasov Poisson equation]{Landau damping for the linearized Vlasov Poisson equation in a weakly collisional regime}

\begin{document}

\author{{\sc Isabelle Tristani}}
\address{Centre de Math\'{e}matiques Laurent Schwartz, \'{E}cole polytechnique, CNRS, Universit\'{e} Paris-Saclay, 91128 Palaiseau Cedex, France.
E-mail: {\tt isabelle.tristani@polytechnique.edu}}



\maketitle

\begin{abstract} 
	In this paper, we consider the linearized Vlasov-Poisson equation around an homogeneous Maxwellian equilibrium in a weakly collisional regime: there is a parameter~$\eps$
	in front of the collision operator which will tend to $0$. Moreover, we study two cases of collision operators, linear Boltzmann and Fokker-Planck. We prove a result of Landau damping for those equations in 
	Sobolev spaces uniformly with respect to the collision parameter $\eps$ as it goes to $0$.  
\end{abstract}

\vspace{0.5cm}
\tableofcontents


\section{Introduction} 
\label{sec:intro}
\setcounter{equation}{0}
\setcounter{theo}{0}


In this paper, we are interested in the phenomenon of Landau damping for a linear Vlasov Poisson equation in a weakly collisional regime. The physical motivation for this problem comes from plasma theory. 
Indeed, in kinetic theory, the Landau equation, which can be derived from the Boltzmann one, allows to describe collisions in plasmas. However, this description is only valid for large times. The Vlasov model, on the contrary,
completely neglects collisions in plasmas. Between these two models, is an intermediate one which is a more realistic physical model and which can be written as follows:
	\beqn \label{eq:general}
		\partial_t f + v \cdot \nabla_x f + F(t,x) \cdot \nabla_v f = \eps \, Q_L(f,f)
	\eeqn
where $f=f(t,x,v)$ is the distribution function in phase space, the position $x$ lies in the torus of size $2\pi$  denoted $\T^d := \R^d/ (2 \pi \Z)^d$, the velocity $v$ in $\R^d$. 
Note that one can reduce to the case where the size of the box is $2 \pi$ by a rescaling argument. 
The operator $Q_L$ is the Landau collision operator (that we do not describe precisely because we won't deal with this general model). The mean-field electrostatic force $F(t,x)$ is given by 
	$$
	F(t,x) = - \nabla W *_x \left(\rho_f(t,x) - \int_{\T^d} \rho_f(t,y) \, dy \right)
	$$
with $\rho_f(t,x) := \int_{\R^d} f(t,x,v) \, dv$. The potential $W(x)$ describes the mean-field interaction between particles, we consider here the case of Coulomb repulsive interactions (which occur between electrons in plasmas) for which
	\beqn \label{eq:assW}
		F = e_0 \, \nabla_x  \Delta_x^{-1} \rho_f \qquad \text{and} \qquad \widehat{W}(k) =  e_0 \, |k|^{-2}, \,\,\, k \in \Z^d \setminus \{0\}
	\eeqn
where $e_0>0$ is the electron charge-to-mass ratio.  

In classical plasma physics, the coefficient $\eps$ equals $(\log \Lambda)/(2 \pi \Lambda)$ 
where $\Lambda$ is the plasma parameter which is very large. It is thus coherent to study this type of model in the limit $\eps \to 0$. We expect the behaviour of solutions to~\eqref{eq:general}
to be complex and depend on the considered scales of time. {When $\eps>0$ is fixed, one can expect the hypocoercive effect of collisions to be dominant over large time $O(1/\eps)$, which should induce an exponential
convergence towards the Maxwellian equilibrium. In the Fokker-Planck case, it would also be interesting to analyse the hypoellipticity of the equation, which would provide a regularization effect. Those phenomena are relevant for fixed $\eps$ and their study is a challenging issue and a very interesting perspective. In the sequel, we are interested only in the phenomenon of Landau damping and in giving a result which is uniform with respect to $\eps$ in the limit $\eps \to 0$ (we recall that when $\eps=0$, Landau damping occurs in time of order $O(1)$).} 

As pointed out by Villani in~\cite{coursVillani2}, the study of this problem is of great interest for several reasons: as we mentioned above, from a physical viewpoint, the model is more realistic than the pure Landau or Vlasov ones; 
it also permits quantification with respect to $\eps$ of the Landau damping; and finally, one could hope to bypass the obstacle of Gevrey regularity that faces the non collision model thanks to regularizing effects induced by collisions {in the case of fixed $\eps$}. 
Let us underline the fact that the last aspect is only significant in the nonlinear case that we won't explore in this paper. 
Our study, which is restricted to the linearized problem with some linear collision operators, can thus be seen as the starting point of a long-range program. 

\smallskip
\subsection{The model}
As mentioned above, we only deal with collision operators denoted $\CC(f)$ which are linear. To simplify the problem, we also suppose that 
	\beqn \label{eq:ker}
	\hbox{Ker} (\CC) = \hbox{Span}(M)
	\eeqn
where $M$ is the Maxwellian distribution 
	\beqn \label{eq:max}
		M(v) :=  (2 \pi)^{-d/2} e^{-|v|^2/2}.
	\eeqn
We consider the model
	\beqn \label{eq:general2}
		\partial_t f + v \cdot \nabla_x f + F(t,x) \cdot \nabla_v f = \eps \, \CC(f).
	\eeqn
We write the solutions of~\eqref{eq:general2} as perturbations of the Maxwellian distribution $M$: 
	$$
	f(t,x,v) = M(v) + h(t,x,v)
	$$
where $h$ is a mean-zero perturbation if $\int_{\T^d \times \R^d} f_{in}(x,v) \, dx \, dv =\int_{\T^d \times \R^d} M(v) \, dx \, dv$, an assumption that we will suppose to hold throughout the paper. Let us emphasize here that we consider a linearization around the Maxwellian distribution $M$ and not around some homogeneous distribution or other as it is usually done for the Vlasov-Poisson equation. It is explained by the fact that, under the assumption~\eqref{eq:ker}, $M$ is also the unique equilibrium of the collision operator $\CC$ satisfying $\int_{\R^d} M(v) \, dv =1$. Then $h$ satisfies the following equation (where we denote $\rho(t,x) := \int_{\R^d} h(t,x,v) \, dv$ the perturbation density):
	\beqn  \label{eq:general3}
		\left\{
			\bal
				&\partial_t h + v \cdot \nabla_x h + F(t,x) \cdot \nabla_v h + F(t,x) \cdot \nabla_v M = \eps \, \CC(h) \\
				&F(t,x) = - \nabla W *_x \rho(t,x) \\
				&h(0,x,v)= h_{in}(x,v).
			\eal
		\right.
	\eeqn
where we have used that $\CC(M)=0$.

In this paper, we only work on the linearized problem associated to~\eqref{eq:general3} which reads:
	\beqn \label{eq:main}
		\left\{
			\bal
				&\partial_t h + v \cdot \nabla_x h + F(t,x) \cdot \nabla_v M = \eps \, \CC(h) \\
				&F(t,x) = - \nabla W *_x \rho(t,x) \\
				&h(0,x,v)= h_{in}(x,v).
			\eal
		\right.
	\eeqn
 We end the description of the model we are concerned with by the definition of the collision operators $\CC$ that we shall consider. It can be a linear Boltzmann operator of type:
	\beqn \label{eq:linearB1}
	\CC(h) := \rho \, M - h
	\eeqn
or a Fokker-Planck operator defined through:
	\beqn \label{eq:FP}
	\CC(h)  := \Delta_v h + \textrm{div}_v (vh).
	\eeqn
These are classical operators of statistical physics that allow to describe the interaction between particles and a fixed background. Note that in the case of the collisional Vlasov equation (without the term coming from the mean-field interaction between particles):
	$$
	\partial_t h + v \cdot \nabla_x h =\eps \, \CC(h)
	$$ 
for fixed $\eps>0$, it is clear that the hypocoercive effect is dominant for large times. We can for example apply the results of Dolbeault, Mouhot and Schmeiser~\cite{DMS} to obtain a result of decay at infinity with a rate of type $e^{-\eps \lambda t}$ for some $\lambda>0$ in some Hilbert space.


\medskip
\subsection{Main results and known results}
\subsubsection*{Main results and comments} 
In this paper, we focus on the limit $\eps \to 0$ and our first objective is to obtain a result of Landau damping in the limit $\eps \to 0$ uniformly with respect to $\eps \in [0,\eps_0]$ for some $\eps_0>0$ for the equation~\eqref{eq:main}.  
Before stating our theorems, we introduce the Sobolev spaces in which we are going to develop our analysis. 

We fix for the rest of the paper $\ell > d/2$ and introduce the space $L^2_x\HH^n_v$, $n \in \N$ associated to the following norm:
	$$
	\|f\|^2_{L^2_x\HH^n_v} := \sum_{|\alpha| \le n} \int_{\T^d \times \R^d} \la v \ra^{2\ell} \, | \partial^\alpha_v f|^2 \, dx \, dv
	$$
where we have used the notation $\la v \ra := (1+|v|^2)^{1/2}$. The dependency of $L^2_x\HH^n_v$ with respect to $\ell$ is not included in the notation since $\ell$ is definitely fixed.  

Here is our main result:

\vspace{-0.15cm}
\begin{theo} 
	Consider a mean-zero distribution $h_{in} \in L^2_x\HH^n_v$ for $n \in \N$, $n \ge 3$. 
	There exist $\eps_0>0$, $\lambda_0>0$ such that for all $\eps \in [0,\eps_0]$, the density $\rho = \rho(t,x)$ of the solution $h=h(t,x,v)$ of~\eqref{eq:main} 
	satisfies the following estimate:
		$$
		{ \forall \, t \ge 0}, \quad \|\rho (t,\cdot) \|_{L^2_x} \le {C \over \la t \ra^{n}},
		$$ 
	the solution $h=h(t,x,v)$ satisfies itself
		$$
		\forall \, t \ge 0, \, \forall \, \xi \in \R^d, \quad |\widehat{h}(t,0,\xi)| \le C \, e^{-\eps \lambda_0t} \\
		$$
	and
		$$
		{ \forall \, t \ge 0}, \,\forall \, k \in \Z^d \setminus \{0\},\, \forall \, \xi \in \R^d, \quad |\widehat{h}(t,k,\xi)| \le {C(\xi) \over \la t \ra^{n-2}}
		$$
	where $C(\xi) := C \, \la \xi \ra^{n-2}$ and $C$ is a positive constant depending on $\|h_{in}\|_{L^2_x \HH^n_v}$. 
\end{theo}

{\begin{rem}
	In the Fokker-Planck case, we have a better estimate which provides us an exponential dissipation (see Theorem~\ref{theo1}) for $\eps$ fixed, 
	which is coherent since it comes from the Fokker-Planck structure, which disappears for $\eps=0$. 
	One could also expect that for $\eps$ fixed, enhanced dissipation occurs over large times with a rate of dissipation of type $e^{-\eps t^3}$ thanks to the diffusion part of the Fokker-Planck operator. 
	However, we are not able to prove such a result, due to the presence of the drift term and the reaction term, the study is more intricate. 
\end{rem}}

%

\begin{rem}
	For $\eps=0$, it is well-known that 
	$$ 
	\forall \, \xi \in \R^d,  \, \forall \, t \ge 0, \quad \widehat{h}(t,0,\xi) = \widehat{h_{in}}(0,\xi),
	$$
	it is not anymore the case for $\eps>0$ due to the presence of the collision operator. 
\end{rem}

Let us underline that our study fits into a Sobolev framework. We do not work with analytic or Gevrey regularity as in~\cite{MV} or~\cite{BMM} in which more regularity is needed to handle the nonlinear problem. However, the study of the linearized problem could have been done in Sobolev regularity, this is what we are doing since we only deal with the linearized problem. 
To do that, we take advantage of the paper by Faou and Rousset~\cite{FR} which consists in the analysis of Landau damping for the (nonlinear) Vlasov-HMF equation in a Sobolev setting. { For this equation, the nonlinear study does not require analytic (or Gevrey regularity),
indeed, due to the fact that the interaction kernel has finite support in Fourier space, the control of the nonlinear effect ``plasma-echo'' is made easier. }
Finally, let us mention that our result can easily be adapted to an analytic framework. Also, for the nonlinear equation~\eqref{eq:general3}, it would be particularly interesting to investigate if
the gain of regularity induced by collisions between particles could allow us to relax the hypothesis on the initial datum. 

Another aspect that we point out is that, to the best of our knowledge, it is the first paper in which such a problematic is investigated for the Vlasov-Poisson equation in a weakly collisional regime (in the limit $\eps \to 0$). 
It should be relevant to compare this kind of question with the one studied by Bedrossian, Masmoudi and Vicol~\cite{BMV} about the two-dimensional Euler equation where the equivalent of $\eps$ should be the inverse of the Reynolds number. 
Concerning this kind of problematic of uniform analysis of large time behaviour with respect to a small parameter, we can also mention the work of the author with Mischler~\cite{MT} in which an analysis of large-time behaviour is developed for several homogeneous Fokker-Planck equations uniformly with respect to a small parameter which allows to go from a model to another. 

\subsubsection*{Strategy of the proof} To deal with this problem, we use the classical strategy to analyse Landau damping. We use method of characteristics and Fourier transform to obtain a closed equation on $\widehat{\rho}(t,k)$ for each $k \in \Z^d$. The presence of the collisional operator raises additional difficulties. Indeed, in this equation, there is a term of the form 
	$$
	\int_0^t K_\eps(t-s) \, \widehat{\rho}(s,k) \, ds
	$$ 
which is typical. However, $K_\eps$ splits into two parts, one coming from the Vlasov-Poisson term and the other coming from the collision term. In order to obtain our final result uniformly with respect to $\eps$, we have to sharply estimate the Fourier transform in time of $K_\eps$ and this is the main technical issue of this paper. We are then able to check that a Penrose type criterion is satisfied uniformly with respect to $\eps$ for $\eps$ small enough. Taking advantage of the analysis developed in~\cite{FR} in a Sobolev framework, we end up with an estimate on $\|{\rho}(t,\cdot)\|_{L^2_x}$. We then go back to the kinetic distribution $h$ and estimate each mode of it separately thanks to the result obtained previously on $\rho$.

\subsubsection*{State of the art}
The behaviour of solutions to the equation~\eqref{eq:main} has been largely studied for some fixed $\eps$. 

If $\eps=0$, we are in the framework of Mouhot and Villani~\cite{MV}, they prove Landau damping (linear and nonlinear) in analytic or Gevrey regularity. 
We also mention the work of Bedrossian, Masmoudi and Mouhot~\cite{BMM} in which a new simplified proof is given in Gevrey regularity. 
Let us also specify that some earlier results were obtained in~\cite{CM} by Caglioti and Maffei and in~\cite{HV} by Hwang and Vel\'{a}zquez. We refer to the works of Ryutov~\cite{Ryu}, Mouhot and Villani~\cite{MV}, Villani~\cite{coursVillani} for a detailed historical background on Landau damping and references therein. 

If $\eps>0$ is fixed, one can expect that the collisional term induces a hypocoercive effect. Some results have been obtained on linear and/or nonlinear frameworks for the Vlasov-Poisson-Fokker-Planck equation. About existence results, we mention
the papers of Victory and O'Dwyer in two dimensions~\cite{VOD}, Rein and Weckler~\cite{RW} in three dimensions for small data and Bouchut~\cite{Bou} in three dimensions. Concerning the long-time behaviour, we point out the papers of Bouchut and Dolbeault~\cite{BD}, Carrillo, Soler and V\'{a}squez~\cite{CSV}, Dolbeaut~\cite{Dol} and H\'{e}rau and Thomman~\cite{HT}, in the latter, both Cauchy problem and long-time behaviour have been studied. 

As far as the case of the linear Boltzmann is concerned, the literature is more scarce and we were not able to find any reference on it. 
However, let us emphasize the fact that a large literature is devoted to the study of the Vlasov-Poisson-Boltzmann equation with a general Boltzmann collision operator. On this topic, we call attention on the paper of Dolbeault and Desvillettes~\cite{DD} which deals with the large time behaviour of solutions. On this matter, we also have to refer to a series of papers of Guo in which a robust energy method is developed to treat the global stability of global Maxwellians in a perturbative framework. There are many variants based on Guo's works and its energy method which allow to analyse the Cauchy problem and the large time behaviour of solutions. Let us mention two papers of Guo on the Vlasov-Poisson-Boltzmann equation~\cite{Guo1,Guo2} the first one falls in a near-vacuum regime, the second one in a near Maxwellian setting. Let us give other references which are concerned with the Cauchy problem in a perturbative setting and/or the large time behaviour of solutions and the rates of convergence, among others, we mention the papers of Duan and Strain~\cite{DS}, Duan, Yang and Zhao~\cite{DYZ}, Duan and Liu~\cite{DL1}, Xia, Xiong and Zhao~\cite{XXZ}. From another viewpoint, an analysis of the spectrum has been carried out in~\cite{LYZ} by Li, Yang and Zhong. There is another phenomenon which could occur when the initial datum of the equation approaches distinct global Maxwellians at far fields, 
the solution to the Cauchy problem does not converge anymore to a constant equilibrium state in large time but we do not enter into details in this regard and refer to the recent paper of Duan and Liu~\cite{DL2} and the references given therein. 

\medskip
\subsection{Notations and tools} \label{subsec:not}
We here define some notations and state an elementary lemma on Fourier transform that we will make good use in what follows. 

\subsubsection*{Fourier transform}
For a function $f=f(x)$, $x \in \T^d$, we define its Fourier transform as follows:
	$$
	\widehat{f} (k) = {1 \over {(2\pi)}^{d/2}} \int_{\T^d} f(x) \, e^{-i x \cdot k} \, dx, \quad k \in \Z^d. 
	$$
Similarly, for a function $f=f(v)$, $v \in \R^d$, we define its Fourier transform by:
	$$
	\widehat{f} (\xi) = {1 \over {(2\pi)}^{d/2}} \int_{\R^d} f(v) \, e^{-i v \cdot \xi} \, dv, \quad \xi \in \R^d. 
	$$
Let us notice that the Maxwellian distribution $M$ defined in~\eqref{eq:max} satisfies $\widehat{M} = M$, an equality that we will use all along the paper. 

\noindent Finally, if $f=f(x,v)$, $(x,v) \in \T^d \times \R^d$, we define its Fourier transform through the following formula:
	$$
	\widehat{f} (k,\xi) = {1 \over {(2\pi)}^{d}} \int_{\R^d} f(x,v) \, e^{-i k \cdot x -i v \cdot \xi} \, dx \, dv, \quad (k,\xi) \in \Z^d \times \R^d. 
	$$
We shall also use the Fourier transform in time, if $f=f(t)$, $t \in \R$, we denote 
	\beqn \label{eq:Fourierintime}
	\widetilde{f}(\tau) = \int_\R f(t) \, e^{-it \tau} \, dt, \quad \tau \in \C. 
	\eeqn


\subsubsection*{Elementary lemma}
We here state a simple result which links Fourier transform and regularity of a function (one can check the proof in~\cite[Lemma~2.1]{FR} in dimension $1$).
\begin{lem} \label{lem:elementary}
	For any $n \in \N$, we have:
		$$
		\forall \,  k \in \Z^d, \quad \forall \, \xi \in \R^d, \quad |\widehat{f}(k,\xi)| \le C(n,\ell) \, \la \xi \ra^{-n} \, \|f\|_{L^2_x\HH^n_v}
		$$
	where $C(n,\ell)$ is a positive constant depending only on $\ell$ and $n$. 
\end{lem}

\subsubsection*{Notation}
Finally, let us specify that we shall use the same notation $C$ for positive constants that may change from line to line, 
when $C$ is a positive constant depending only on fixed number.

\medskip
\subsection{Outline of the paper}
In Section~\ref{sec:linearB1}, we deal with the behaviour of solutions of~\eqref{eq:main} with the linear Boltzmann operator~\eqref{eq:linearB1} as a collision operator. In Section~\ref{sec:FP}, we handle the problem with collisions given by a Fokker-Planck operator~\eqref{eq:FP}. 

\medskip
\noindent\textbf{Acknowledgments.} The author would like to thank Daniel Han-Kwan for enlightened discussions, his help and his suggestions. 
The author has been supported by the {\em Fondation Math\'{e}matique Jacques Hadamard}.  

\bigskip

\section{Linear Boltzmann type collisions} \label{sec:linearB1}
\setcounter{equation}{0}
\setcounter{theo}{0}

In this section, we deal with the following model: 
	\beqn \label{B1eq:VPB1}
		\left\{
			\bal
				&\partial_t h + v \cdot \nabla_x h + F(t,x) \cdot \nabla_v M = \eps \, (\rho \, M - h) \\
				&F(t,x) = - \nabla W *_x \rho(t,x) \\
				&h(0,x,v)= h_{in}(x,v).
			\eal
		\right.
	\eeqn
where we recall that $\rho(t,x) = \int_{\R^d} h(t,x,v) \, dv$.
The main theorem reads:
\vspace{-0.15cm}
\begin{theo} \label{theo0}
	Consider a mean-zero distribution $h_{in} \in L^2_x\HH^n_v$ for $n \in \N$, $n \ge 3$. 
	There exists $\eps_0>0$ such that the solution $h=h(t,x,v)$ of \eqref{B1eq:VPB1} and its associated density $\rho = \rho(t,x)$ satisfy the following estimates:
		\beqn \label{B1eq:finalrho}
			\forall \, \eps \in [0,\eps_0], \, \forall \, t \ge 0, \quad \|\rho (t,\cdot) \|_{L^2_x} \le {C \over \la t \ra^{n}},
		\eeqn 
	and
		\beqn \label{B1eq:finalh}
			\left\{
			\bal
				&\forall \, \eps \ge 0, \, \forall \, \xi \in \R^d, \, \forall \, t \ge 0, \quad |\widehat{h}(t,0,\xi)| \le C \, e^{-\eps t}, \\
				&\forall \, \eps \in [0,\eps_0], \,  \forall \, k \in \Z^d \setminus \{0\},\, \forall \, \xi \in \R^d, \, \forall \, t \ge 0, \quad |\widehat{h}(t,k,\xi)| \le {C(\xi) \over \la t \ra^{n-2}},
			\eal
			\right.
		\eeqn
	where $C(\xi) := C \, \la \xi\ra^{n-2}$ and $C$ is a positive constant depending on $\|h_{in}\|_{L^2_x\HH^n_v}$.
\end{theo}

\medskip
\subsection{The Penrose criterion}
Before going into the proof of Theorem~\ref{theo0}, we start by giving a lemma which allows us to apply a result from~\cite{FR}. To do that, we introduce the following notations: for $(t,k) \in \R \times \Z^d$, we define for $\eps \ge 0$,
	\beqn \label{B1def:Keps}
	\left\{
		\bal
			&K^0_\eps(t,k) := \eps \, e^{-\eps t} \, \widehat{M}(kt) \, \mathds{1}_{t \ge 0}\\
			&K^1_\eps(t,k) := - \widehat{W}(k)\, e^{-\eps t}  \, \widehat{M}(kt) \, |k|^2 \, t \, \mathds{1}_{t \ge 0} \\
			&K_\eps (t,k) := K^0_\eps(t,k) + K^1_\eps(t,k).
		\eal
	\right. 
	\eeqn
	
\begin{lem} \label{B1lem:Keps}
	There exists $\eps_0>0$ such that the kernel $K_\eps$ defined in~\eqref{B1def:Keps} satisfies the following condition:		
		$$
		\emph{\bf{(H)}} \qquad \qquad \exists \, \kappa>0, \quad \forall \, \eps \in [0,\eps_0], \quad  \inf_{k \in \Z^d \setminus \{0\}} \inf_{\emph{Im} \, \tau \le 0} |1- \widetilde{K_\eps}(\tau,k)| \ge \kappa
		$$
	where we recall that $\widetilde{K_\eps}(\tau,k)$ stands for the Fourier transform in time of $K_\eps$ defined in~\eqref{eq:Fourierintime}.
\end{lem}

\noindent {\it Proof of Lemma~\ref{B1lem:Keps}.}
	We start by giving an estimate on $|1- \widetilde{K_0^1}(\tau,k)|$. We recall that since all marginals of $M$ are increasing/decreasing, 
	one can check (see~\cite{coursVillani}) that the gaussian distribution $M$ satisfies the following form of Penrose stability criterion~\cite{Penrose} (given in \cite{MV} in this form):
		$$
		\forall \, k \in \Z^d \setminus \{0\}, \quad \forall \, w \in \R \, \, \text{s.t.} \, \, (M_k)'(w) =0, \quad \widehat{W}(k)  \, \left( \hbox{p.v.} \int_\R \frac{(M_k)'(r)}{r- w}\, dr \right)  <1
		$$
	where $M_k$ denotes the marginals of $M$ along $k \in \Z^d\setminus\{0\}$:
		$$
		M_k(r) := \int_{kr/ |k| + k^\perp} M(w) \, dw.
		$$
	Consequently, we know (see~\cite{MV} and~\cite{BMM} for a complete proof) that there exists $\lambda_\dagger>0$ and $\kappa_0>0$ such that 
		\beqn \label{B1eq:conditionL}
			\forall \, \xi \in \C, \,  \, \Re e \, \xi < \lambda_\dagger, \quad \inf_{k \in \Z^d} |\LL(k,\xi) - 1| \ge \kappa_0
		\eeqn
	where $\LL$ is defined through ($\bar \xi$ is the complex conjugate of $\xi$)
		$$
		\LL(k,\xi) = - \int_0^\infty e^{\bar \xi |k| t} \widehat{M}(kt) \, \widehat{W}(k) \, |k|^2 \, t \, dt, \quad (k,\xi) \in \Z^d \times \R^d,
		$$
	note that this integral is well defined because it is absolutely convergent due to the decay property of $\widehat{M} = M$. 
	Then, one can check that if we denote $\tau = \lambda + i \zeta$ with $(\lambda, \zeta) \in \R^2$, 
		$$
		\widetilde{K_0^1}(\tau,k) = \LL\left( k,\frac{\zeta+i\lambda}{|k|}\right)
		$$
	and thus using~\eqref{B1eq:conditionL}, we deduce that for any $k \in \Z^d \setminus \{0\}$, 
		\beqn \label{B1eq:K0}
		\forall \, \tau \in \C, \, \, \hbox{Im} \tau \le 0, \quad  |\widetilde{K_0^1}(\tau,k) - 1| \ge \kappa_0.
		\eeqn
		
	Still using the notation $\tau = \lambda + i \zeta$, we also have 
		$$
		|\widetilde{K^0_\eps}(\tau,k)| =  \eps \, \left|\int_0^\infty e^{-it\tau}  \, e^{-\eps t} \, \widehat{M}(kt) \, dt\right| \le \eps \, \int_0^\infty e^{\zeta t}  \, e^{-\eps t} \, M(kt) \, dt.
		$$
	We deduce for any $k \in \Z^d \setminus \{0\}$, 
		\beqn \label{B1eq:K0eps}
			\forall \, \tau \in \C, \, \, \hbox{Im} \tau \le 0, \quad |\widetilde{K^0_\eps}(\tau,k)| \le \eps \, \int_0^\infty M(kt) \, dt \le {\eps \over |k|} \, \int_0^\infty M(t) \, dt \le C \, \eps.
		\eeqn
	
	Then, it remains to estimate the difference coming from $K^1_\eps - K^1_0$. We have:
		$$
		\bal
			|\widetilde{K^1_\eps}(\tau,k) - \widetilde{K^1_0}(\tau, k)| &= \left| \int_0^\infty e^{-it\tau} \, (K^1_\eps(t,k) - K^1_0(t,k)) \, dt \right| \\
			&=  e_0 \, \left| \int_0^\infty e^{-it\tau} \, (1-e^{-\eps t}) \, \widehat{M}(kt) \, t \, dt\right| \\
			&\le C \,  \eps \, \int_0^\infty e^{\zeta t}  \, t^2 \, e^{-\frac{|k|^2 t^2}{2}} \, dt
		\eal
		$$
	where we used that $1-e^{-\eps t} \le \eps \, t$ for any $t \ge 0$.
	Performing a change of variable in the last integral, we get for any $k \in \Z^d \setminus \{0\}$,
		\beqn \label{B1eq:K1diff}
			\forall \,  \tau \in \C, \, \, \hbox{Im} \tau \le 0, \quad |\widetilde{K^1_\eps}(\tau,k) - \widetilde{K^1_0}(k, \tau)| \le \frac{C}{|k|^3} \, \eps \le C \, \eps.
		\eeqn
		
	We now go back to our full problem which consists in estimating the difference 
		$$
		|1-\widetilde{K_\eps}(\tau,k)| \quad \text{with} \quad K_\eps = K^0_\eps + K^1_0 + K^1_\eps - K^1_0.
		$$
	Gathering~\eqref{B1eq:K0},~\eqref{B1eq:K0eps} and~\eqref{B1eq:K1diff}, we conclude that for any $k \in \Z^d \setminus \{0\}$, for any $\tau \in \C$, $\hbox{Im} \tau \le 0$,
		$$
		\bal
			|\widetilde{K_\eps}(\tau,k) - 1| &\ge |\widetilde{K_0^1}(\tau,k) - 1| - |\widetilde{K^1_\eps}(\tau,k) - \widetilde{K^1_0}(k, \tau)| - |\widetilde{K^0_\eps}(\tau,k)| \\
			&\ge \kappa_0 - c_0 \eps,
		\eal
		$$
	for some positive constant $c_0>0$ which only depends on a fixed number.
	We then choose $\eps_0>0$ such that $\kappa_0 - c_0 \eps_0 >0$. Denoting $\kappa := \kappa_0 - c_0 \eps_0$, we obtain
		$$
		 \forall \, \eps \in [0,\eps_0], \quad \inf_{k \in \Z^d \setminus \{0\}} \inf_{\hbox{Im} \, \tau \le 0} |\widetilde{K_\eps}(\tau,k) - 1| \ge \kappa,
		$$
	from which we conclude that {\bf(H)} is fulfilled. 
\qed

\medskip
\subsection{Proof of Theorem~\ref{theo0}}
We can now handle the proof of Theorem~\ref{theo0}, we divide it into two parts, the first one being devoted to the proof of estimate~\eqref{B1eq:finalrho} and the second one to~\eqref{B1eq:finalh}.

\smallskip
\noindent {\it Proof of the estimate~\eqref{B1eq:finalrho}.} 
	The first thing that we have to notice is that, due to the conservation properties of the operator $- v \cdot \nabla_x h - F(t,x) \cdot \nabla_v M + \eps \, (\rho \, M - h)$, we know that for all times, 
	$\widehat{\rho}(t,0) = (2 \pi)^{d/2} \, \widehat{h_{in}}(0,0)=0$ since $h_{in}$ is a mean-zero distribution. 
	
	We then focus on the study of $\widehat{\rho}(t,k)$ for $k \in \Z^d \setminus \{0\}$. In the rest of the proof, we consider $\eps \in [0,\eps_0]$ where $\eps_0$ is given in Lemma~\ref{B1lem:Keps}. 
	Using the method of characteristics, we have 
		$$
		\bal
			h(t,x,v) &= h_{in}(x-vt,v) \, e^{-\eps t} + \eps \, \int_0^t e^{-\eps (t-s)} \, \rho(s,x-v(t-s)) \, M(v) \, ds \\
			&\quad - \int_0^t e^{-\eps(t-s)} \, (\nabla_v M \cdot F)(s,x-v(t-s),v) \, ds \\
		\eal
		$$
	We then take the Fourier transform in both variables $(x,v)$ of the previous equality:
		$$
		\bal
			&\quad \widehat{h}(t,k,\xi) = e^{-\eps t} \, \int_{\T^d \times \R^d} h_{in}(x-vt,v) \, e^{-i k \cdot x -i v \cdot \xi} \, dx \, dv \\
			& + {\eps \over (2\pi)^d} \int_0^t \int_{\T^d \times \R^d} e^{-\eps(t-s)} \,\rho(s,x-v(t-s)) \, M(v) \, e^{-i k \cdot x -i v \cdot \xi} \, dx \, dv \, ds \\
			& - {1 \over (2\pi)^d} \int_0^t \int_{\T^d \times \R^d} e^{-\eps(t-s)} \, (\nabla_v M \cdot F)(s,x-v(t-s),v) \, e^{-i k \cdot x -i v \cdot \xi} \, dx \, dv \, ds.
		\eal
		$$
	Performing change of variables, we obtain
		$$
		\bal
			&\quad \widehat{h}(t,k,\xi) = e^{-\eps t} \,  \widehat{h_{in}} (k, \xi + kt) \\
			& + {\eps \over (2\pi)^d} \int_0^t \int_{\T^d \times \R^d} e^{-\eps(t-s)} \,\rho(s,x) \, M(v) \, e^{-i k \cdot x -i v \cdot (\xi + k(t-s))} \, dx \, dv \, ds \\
			&- {1 \over (2\pi)^d} \int_0^t \int_{\T^d \times \R^d} e^{-\eps(t-s)} \,  (\nabla_v M \cdot F)(s,x,v) \, e^{-i k \cdot x -i v \cdot (\xi + k(t-s))} \, dx \, dv \, ds,
		\eal
		$$
	from which we deduce
		\beqn \label{B1eq:Fourier1}
		\bal
			\widehat{h}(t,k,\xi) &= e^{-\eps t}  \, \widehat{h_{in}} (k, \xi + kt) + \eps \int_0^t e^{-\eps(t-s)} \, \widehat{\rho}(s,k) \, \widehat{M}(\xi + k(t-s)) \, ds \\
			&\quad - \int_0^t e^{-\eps(t-s)} k \cdot (\xi + k(t-s)) \, \widehat{W}(k) \, \widehat{\rho}(s,k) \, \widehat{M}(\xi+ k(t-s)) \, ds.
		\eal
		\eeqn
	Taking $\xi=0$ in~\eqref{B1eq:Fourier1}, we obtain the closed equation on $\widehat{\rho}(t,k)$:
		$$
		\bal
			\widehat{\rho}(t,k) &= e^{-\eps t}  \, \widehat{h_{in}} (k, kt) + \eps \int_0^t e^{-\eps(t-s)} \, \widehat{\rho}(s,k) \, \widehat{M}(k(t-s)) \, ds \\
			&\quad - \int_0^t e^{-\eps(t-s)} k \cdot ( k(t-s)) \, \widehat{W}(k) \, \widehat{\rho}(s,k) \, \widehat{M}(k(t-s))  \, ds.
		\eal
		$$
	that we can sum up in
		$$
			\widehat{\rho}(t,k) = e^{-\eps t}  \, \widehat{h_{in}} (k, kt) + \int_0^t K_\eps (t-s) \,  \widehat{\rho} (s,k) \, ds 
		$$
	where we recall that $K_\eps$ is defined in~\eqref{B1def:Keps}. 
	We can notice here that we are almost in the same framework as the one of \cite[Lemma~2.3]{FR} for $\eps\in [0,\eps_0]$. Indeed, we have $\la \cdot \ra^2 M \in L^2_x\HH^n_v$ for any $n \in \N$ 
	and from Lemma~\ref{B1lem:Keps}, we know that $K_\eps$ satisfies the following condition:		
		$$
		\textbf{(H)} \qquad \qquad \exists \, \kappa>0, \quad \forall \, \eps \in [0,\eps_0], \quad  \inf_{k \in \Z^d \setminus \{0\}} \inf_{\hbox{Im} \, \tau \le 0} |1- \widetilde{K_\eps}(\tau,k)| \ge \kappa.
		$$
	The only difference is that the proof in~\cite{FR} is done only for the case $k=\pm 1$. But it is uncomplicated to extend the result for $|k|>1$. We then obtain that for any $\gamma \ge 0$, for any $T \ge 0$:
		\beqn \label{B1eq:FR}
			\forall \, k \in \Z^d \setminus \{0\}, \quad \sup_{t \in [0,T]} \la kt \ra^\gamma \, |\widehat{\rho}(t,k)| \le \sup_{t \in [0,T]} \la kt \ra^\gamma \, e^{-\eps t}  \, |\widehat{h_{in}} (k, kt)|  
			\le \sup_{t \in [0,T]} \la kt \ra^\gamma\, |\widehat{h_{in}} (k, kt)|.
		\eeqn
	Using Lemma~\ref{lem:elementary}, there holds,
		\beqn \label{B1eq:FR2}
			|\widehat{h_{in}} (k, kt)| \le C(n) \,  \la kt \ra^{-n} \, \|h_{in}\|_{L^2_x\HH^n_v}.
		\eeqn
	Combining~\eqref{B1eq:FR} and~\eqref{B1eq:FR2} with $\gamma=n$, we obtain for any $T \ge 0$,
		$$
		\forall \, k \in \Z^d \setminus \{0\}, \quad \sup_{t \in [0,T]} \la kt \ra^n \, |\widehat{\rho}(t,k)| \le  C(n) \, \|h_{in}\|_{L^2_x\HH^n_v}.
		$$
	from which we deduce		
		\beqn \label{B1eq:rho}
			\forall \, k \in \Z^d \setminus \{0\}, \, \forall \, t \ge 0, \quad |\widehat{\rho}(t,k)| \le  {C(n) \over \la kt \ra^n} \, \|h_{in}\|_{L^2_x\HH^n_v},
		\eeqn
	this yields the final result~\eqref{B1eq:finalrho}. 
\qed

\smallskip
We now concentrate on the proof of estimate on the solution~$h=h(t,x,v)$.

\smallskip
\noindent {\it Proof of the estimate~\eqref{B1eq:finalh}.} 
	The first part concerns the moment of order $0$ in $x$. Contrary to the case of the linear Vlasov-Poisson equation, the latter is not preserved due to the presence of the collision operator. 
	Indeed, noticing that since $h$ is mean-zero, we have:
		$$
		\partial_t \widehat{h}(t,0,\xi) = -\eps \, \widehat{h}(t,0,\xi)
		$$
	we then deduce that for any fixed $\xi \in \R^d$,
		$$
		\forall \, \eps \ge 0, \quad |\widehat{h}(t,0,\xi)| = e^{-\eps t} \, |\widehat{h_{in}}(0,\xi)| \le C \, e^{-\eps t} \, \|h_{in}\|_{L^2_x \HH^n_v} 
		$$
	where the last estimate comes from Lemma~\ref{lem:elementary}.

	We now deal with the case $k \in \Z^d \setminus \{0\}$ and consider $\eps \in[0,\eps_0]$. From~\eqref{B1eq:Fourier1}, we have
		$$
		\bal
			\widehat{h}(t,k,\xi-kt) &= e^{-\eps t}  \, \widehat{h_{in}} (k, \xi) + \eps \int_0^t e^{-\eps(t-s)} \, \widehat{\rho}(s,k) \, \widehat{M}(\xi -ks) \, ds \\
			&\quad - \int_0^t e^{-\eps(t-s)} k \cdot (\xi -ks) \, \widehat{W}(k) \, \widehat{\rho}(s,k) \, \widehat{M}(\xi-ks) \, ds \\
			&=: T_1+T_2+T_3. 
		\eal
		$$
	Since $h_{in} \in L^2_x\HH^n_v$, we have from Lemma~\ref{lem:elementary} that 
		\beqn \label{B1eq:T1}
			\forall \, k \in \Z^d \setminus \{0\}, \, \forall \, \xi \in \R^d, \quad |T_1| \le C \, \la\xi\ra^{-n}. 
		\eeqn 
	To estimate $T_2$ and $T_3$, we use the estimate~\eqref{B1eq:rho} on $\rho$ that we obtained previously combined with the decay property of $\widehat{M}=M$ which implies for example
		$$
		 \widehat{M}(\xi - ks) \le C \, \la \xi - ks \ra^{2-n}
		$$
	Concerning $T_2$, we thus have
		$$
		|T_2| \le C  \int_0^t \la ks \ra^{-n} \, \la \xi - ks \ra^{2-n} \, ds \le C  \int_0^t \la ks \ra^{-n} \, \la \xi \ra^{2-n} \, \la ks \ra^{n-2} \, ds
		$$
	where we have used the estimate $\la \xi - ks \ra^{2-n} \le C \, \la \xi \ra^{2-n} \, \la ks \ra^{n-2}$, we deduce
		$$
		|T_2| \le C \, \la \xi \ra^{2-n}  \int_0^t  \frac{ds}{\la |k| s \ra^2}.
		$$
	Performing a change of variable, we obtain 
		\beqn \label{B1eq:T2}
			\forall \, k \in \Z^d \setminus \{0\}, \quad \forall \, \xi \in \R^d, \quad |T_2|  \le C \, {1 \over |k|} \, \la \xi \ra^{2-n} \le C \, \la \xi \ra^{2-n}.
		\eeqn
	We treat $T_3$ in a similar way, noticing that $|k| \, \widehat{W}(k) = O(1)$, and that we have 
		$$
		|\xi - ks| \, \widehat{M}(\xi - ks) \le C \, \la \xi - ks \ra^{2-n}.
		$$
	We conclude this proof as for $T_2$ and we obtain 
		\beqn \label{B1eq:T3}
			\forall \, k \in \Z^d \setminus \{0\}, \, \forall \, \xi \in \R^d,  \, \forall \, t \ge0,\quad |T_3|  \le  C \, \la \xi \ra^{2-n}.
		\eeqn
	Combining~\eqref{B1eq:T1},~\eqref{B1eq:T2} and~\eqref{B1eq:T3}, we obtain 
		$$
		\forall \, k \in \Z^d \setminus \{0\}, \, \forall \, \xi \in \R^d, \, \forall \, t \ge0, \quad |\widehat{h}(t,k,\xi-kt)| \le {C \over \la\xi\ra^{n-2}}. 
		$$
	In particular, we have for any fixed $\xi \in \R^d$,
		$$
		\forall \, k \in \Z^d \setminus \{0\},  \, \forall \, t \ge0,\quad |\widehat{h}(t,k,\xi)| \le {C \over \la\xi+kt\ra^{n-2}} \le {C(\xi) \over \la t \ra^{n-2}}
		$$
	where we have used that $\la\xi+kt\ra^{2-n} \le C \,\la \xi \ra^{n-2} \, \la kt \ra^{2-n} \le C \, \la \xi \ra^{n-2} \, \la t \ra^{2-n}$
	and this concludes the proof of~\eqref{B1eq:finalh}. 
\qed


\bigskip

\section{Fokker-Planck type collisions} \label{sec:FP}
\setcounter{equation}{0}
\setcounter{theo}{0}

In this section, we deal with the following model: 
	\beqn \label{FPeq:VPFP}
		\left\{
			\bal
				&\partial_t h + v \cdot \nabla_x h + F(t,x) \cdot \nabla_v M = \eps \, (\Delta_v h + \hbox{div}_v(vh)) \\
				&F(t,x) = - \nabla W *_x \rho(t,x) \\
				&h(0,x,v)= h_{in}(x,v)
			\eal
		\right.
	\eeqn	
where we recall that $\rho(t,x) = \int_{\R^d} h(t,x,v) \, dv$.

\begin{theo}[Fokker-Planck operator] \label{theo1}
	Consider a mean-zero distribution $h_{in} \in L^2_x\HH^n_v$ for $n \in \N$, $n \ge 3$. 
	There exists $\eps_0>0$ such that the solution $h=h(t,x,v)$ of~\eqref{FPeq:VPFP}
	and its associated density $\rho = \rho(t,x)$ satisfy the following estimates:
		\beqn \label{FPeq:finalrho}
		\bal
			&\forall \, \eps \in [0,\eps_0], \, { \forall \, t \ge0}, \quad { \|\rho (t,\cdot) \|_{L^2_x} \le C \, {e^{-\eps t} \over \langle t \rangle^n} \le {C \over \langle t \rangle^n}},
		\eal
		\eeqn
	and
		\beqn \label{FPeq:finalh}
			\left\{
			\bal
				&\forall \, \eps \ge0, \, \forall \, \xi \in \R^d, \, \forall \, t \ge 0, \quad |\widehat{h}(t,0,\xi)| \le  C \, e^{-\eps \lambda_* t}, \\
				&\forall \, \eps \in [0,\eps_0], \,  \forall \, k \in \Z^d \setminus \{0\},\, \forall \, \xi \in \R^d, \,  { \forall \, t \ge 0}, \quad { |\widehat{h}(t,k,\xi)|  
				\le C(\xi){e^{-\eps t} \over \la t \ra^{n-2}} \le {C(\xi) \over \la t \ra^{n-2}}}, \\
			\eal
			\right.
		\eeqn
	where $\lambda_*>0$ is the rate of decay for the homogeneous Fokker-Planck equation  in $L^2(\la v \ra^\ell)$ (see Lemma~\ref{FPlem:homFP}) and $C(\xi) := C \, \la \xi\ra^{n-2}$ 
	and $C$ is a positive constant depending on $\|h_{in}\|_{L^2_x\HH^n_v}$. 
	

\end{theo}

\medskip
\subsection{The Penrose criterion}
As in Section~\ref{sec:linearB1}, we introduce a kernel $K_\eps$ which is going to be useful in the proof of Theorem~\ref{theo1}. For $(t,k) \in \R \times \Z^d$ and $\eps \in (0,1/12]$, we define
	{ 
	\beqn \label{FPdef:Keps}
			K_\eps (t,k):= - e^{\eps t} \, \hbox{exp}\left(-\eps \int_0^t\left|\chi_\eps(\sigma) \, k \right|^2 \, d\sigma \right)  \, k \cdot \left(\chi_\eps(t) \, k\right) \, 
				\widehat{W}(k) \, \widehat{M}\left(\chi_\eps(t) \, k\right) \, \mathds{1}_{t \ge 0}, 
	\eeqn}
where we denoted for $\eps>0$,
	$$
	\chi_\eps(t) := \frac{1-e^{-\eps t}}{\eps}, \quad t \in \R.
	$$
We also introduce a kernel $K_0$ to which we shall compare $K_\eps$:
	\beqn \label{FPdef:K0}
		K_0 (t,k) := - \widehat{W}(k) \, |k|^2 \, M(kt) \, t \, \mathds{1}_{t \ge 0}. 	
	\eeqn
The aim of the next part is to prove that $K_\eps$ satisfies a Penrose type criterion such as we did in Section~\ref{sec:linearB1}, Lemma~\ref{B1lem:Keps}. 

\begin{lem} \label{FPlem:Keps}
	There exists $\eps_0>0$ such that the kernel $K_\eps$ defined in~\eqref{FPdef:Keps} satisfies the following condition:		
		$$
		\emph{\bf{(H)}} \qquad \qquad \exists \, \kappa>0, \quad \forall \, \eps \in [0,\eps_0], \quad  \inf_{k \in \Z^d \setminus \{0\}} \inf_{\emph{Im} \, \tau \le 0} |1- \widetilde{K_\eps}(\tau,k)| \ge \kappa
		$$
	where we recall that $\widetilde{K_\eps}(\tau,k)$ stands for the Fourier transform in time of $K_\eps$ defined in~\eqref{eq:Fourierintime}.
\end{lem}

\noindent {\it Proof of Lemma~\ref{FPlem:Keps}.} First, we recall that we have shown in the proof of Lemma~\ref{B1lem:Keps} (see estimate~\eqref{B1eq:K0}) that for any $k \in \Z^d \setminus \{0\}$, 
		\beqn \label{FPeq:K0}
			\forall \, \tau \in \C, \, \, \hbox{Im} \tau \le 0, \quad  |\widetilde{K_0}(\tau,k) - 1| \ge \kappa_0
		\eeqn
	for some $\kappa_0>0$.
	
	We then estimate the difference $\widetilde{K_\eps}(/tau,k) - \widetilde{K_0}(\tau,k)$ for $k \in \Z^d \setminus \{0\}$ splitting it into several parts. 
	For any $\tau = \lambda + i \zeta$, $\hbox{Im} \tau = \zeta \le 0$, we have:
		{ $$
		\bal
			&\quad| \widetilde{K}_\eps(\tau,k) - \widetilde{K}_0(\tau,k)| \\
			&\le e_0 \int_0^\infty \left| \chi_\eps (t) \,M(k\chi_\eps(t)) - t \, M(kt) \right| \, \mathds{1}_{0 \le t \le 3/\eps} \, e^{\eps t} \, 
			\hbox{exp} \left( - \eps \int_0^t \chi^2_\eps(\sigma) \, |k|^2 \, d\sigma \right) \, |e^{-it\tau}| \, dt \\
			&\quad + e_0 \int_0^\infty \left| e^{\eps t} \,  \hbox{exp} \left( - \eps \int_0^t \chi^2_\eps(\sigma) \, |k|^2 \, d\sigma \right) -1 \right| t \, M(kt) \, \mathds{1}_{0 \le t \le 3/\eps} \, |e^{-it\tau}|  \, dt \\
			&\quad + e_0 \int_0^\infty e^{\eps t} \,  \hbox{exp} \left( - \eps \int_0^t \chi^2_\eps(\sigma) \, |k|^2 \, d\sigma \right)  \, \left| 
			\chi_\eps (t) \,M(k\chi_\eps(t)) - t \, M(kt) \right|\, \mathds{1}_{ t \ge 3/\eps} \, |e^{-it\tau}| \, dt \\
			&\quad + e_0  \int_0^\infty \left| e^{\eps t} \, \hbox{exp} \left( - \eps \int_0^t \chi^2_\eps(\sigma) \, |k|^2 \, d\sigma \right) - 1 \right| \, t\, M(kt) \, \mathds{1}_{ t \ge 3/\eps} \, |e^{-it\tau}| \, dt \\
			&=: I_1 + I_2 + I_3+ I_4. 
		\eal
		$$}
	
	{ Concerning the term $I_1$, we have since $\zeta = \hbox{Im}\,  \tau \le 0$ and $e^{\eps t} \le e^3$ for $t \le 3/\eps$, 
		$$
		\bal
			I_1 &\le {C } \, \int_0^{3/\eps} \left|  \chi_\eps (t) \,M(k\chi_\eps(t)) -  t \, M(kt) \right|   \,e^{\zeta t} \, dt \\
			&\le C \, \int_0^{3/\eps} \left|  \chi_\eps (t)  -  t \right| \, M(kt)  \, dt + {C} \, \int_0^{3/\eps} \left| M(k\chi_\eps(t)) -  M(kt) \right|  \, \chi_\eps(t) \, dt \\
			&=: I_{11} + I_{12}. 
		\eal
		$$ 
	First, we perform a straightforward analysis of the function $t \mapsto t-\chi_\eps(t) - \eps t^2/2$ to prove that 
		\beqn \label{FPeq:chieps}
			\forall \, t \in \R^+, \quad |\chi_\eps(t)-t| = t - \chi_\eps(t) \le \eps\, {t^2 \over 2}. 
		\eeqn
	We deduce that 
		\beqn \label{FPeq:I11}
		\bal
			I_{11} &\le C \, \eps \,  \int_0^{3/\eps} t^2 \, M(kt)   \, dt \le { {C \, \eps} \over |k|^2} \,  \int_0^{3/\eps} |k|^2 \, t^2 \, M(kt)   \, dt  \\
			&\le {{C \, \eps}\over |k|^3} \int_0^{3|k|/\eps} t^2 \, e^{-t^2/2} \, dt \le C \, \eps \int_0^{\infty} t^2 \, e^{-t^2/2} \, dt \le C \, \eps
		\eal
		\eeqn
	where we performed a change of variable. 
	Concerning $I_{12}$, we notice that 
		$$
		\left| M(k\chi_\eps(t)) -  M(kt) \right| = \left| M^{1/2}(k\chi_\eps(t)) -  M^{1/2}(kt) \right| \, \left(M^{1/2}(k\chi_\eps(t)) +  M^{1/2}(kt) \right).
		$$
	But, one can easily prove that $\chi_\eps(t) \ge e^{-3}  t$ for any $t \in [0,3/\eps]$. We thus have, denoting $\mu := e^{- {{|\cdot|^2} \over {4 e^{6}}}}$,
		$$
		M^{1/2}(k\chi_\eps(t)) +  M^{1/2}(kt) \le C  \, \mu(kt).
		$$
	Also, the gradient of $M^{1/2}$ is bounded in $\R^d$, we deduce that for any $t \in [0,3/\eps]$,
		$$
		\bal
			&\left| M(k\chi_\eps(t)) -  M(kt) \right| \le C \, \left| M^{1/2}(k\chi_\eps(t)) -  M^{1/2}(kt) \right| \, \mu(kt) \\
			&\qquad \quad \le C \, |k| (t - \chi_\eps (t)) \, \mu(kt) \le C \, \eps \,|k| \, t^2 \, \mu(kt)
		\eal
		$$
	where the last inequality comes from~\eqref{FPeq:chieps}. 
	This implies, since $\chi_\eps(t) \le t$ for any $t \ge 0$:
		\beqn \label{FPeq:I12}
		\bal
			I_{12} &\le {C} \, \eps \,\int_0^{3/\eps} |k| \, t^3 \,  \mu(kt)  \, dt \le {{C \, \eps} \over |k|^3} \int_0^{6|k| / {\eps}} t^3 \, e^{-\left({t \over {2e^6}}\right)^2} \, dt \\
			&\le {{C \, \eps} \over |k|^3} \int_0^{\infty} t^3 \, e^{-\left({t \over {2e^6}}\right)^2} \, dt \le C \, \eps.
		\eal
		\eeqn
	}
	{ Concerning $I_2$, we also split its study into two parts: 
		$$
			\bal
				I_2 &\le C \,  \int_0^\infty  \left|  \hbox{exp} \left( - \eps \int_0^t \chi^2_\eps(\sigma) \, |k|^2 \, d\sigma \right) -1 \right| t \, M(kt) \, \mathds{1}_{0 \le t \le 3/\eps} \,e^{\zeta t} \,  dt \\
				&\quad + C \,  \int_0^\infty e^{\eps t} \, (1- e^{-\eps t}) \,   \hbox{exp} \left( - \eps \int_0^t \chi^2_\eps(\sigma) \, |k|^2 \, d\sigma \right) t \, M(kt) \, \mathds{1}_{0 \le t \le 3/\eps}\, e^{\zeta t} \, dt \\
				&=: I_{21} + I_{22}. 
			\eal
		$$
	To deal with $I_{21}$, we study the function $t \mapsto 1-  \hbox{exp} \left( - \eps \int_0^t \chi^2_\eps(\sigma) \, |k|^2 \, d\sigma \right) - \eps \, |k|^2 \, t^3/3$ to obtain
		$$
		\bal
			\forall \, t \ge 0, &\quad \left|  \hbox{exp} \left( - \eps \int_0^t \chi^2_\eps(\sigma) \, |k|^2 \, d\sigma \right) -1 \right|\\
			&= 1 - \hbox{exp} \left( - \eps \int_0^t \chi^2_\eps(\sigma) \, |k|^2 \, d\sigma \right) \le  \eps \, |k|^2 \, {t^3 \over 3}.
		\eal
		$$
	We deduce that for any $k \in \Z^d \setminus \{0\}$,
		\beqn \label{FPeq:I21}
			I_{21} \le C \, \eps \, |k|^2 \, \int_0^\infty t^4 \, M(kt) \, e^{\zeta t} \, dt \le {{C \, \eps} \over |k|^3}  \int_0^\infty t^4 \, e^{-t^2/2} \, dt \le C \, \eps,
		\eeqn
	where we have performed a change of variable to get the last inequality but one. For $I_{22}$, we use the simple inequality 
		$$
		\forall \, t \ge 0, \quad1-e^{-\eps t} \le \eps \, t
		$$
	so that 
		\beqn \label{Fpeq:I22}
			I_{22} \le C \, \eps \, \int_0^\infty e^{3} \, t^2 \, M(kt) \, e^{\zeta t}\, dt \le {{C \, \eps} \over |k|^3}  \int_0^\infty t^2 \, e^{-t^2/2} \, dt \le C \, \eps. 
		\eeqn
	}
	{ We then control $I_3$ exploiting the term $\hbox{exp} \left( - \eps \int_0^t \chi^2_\eps(\sigma) \, |k|^2 \, d\sigma \right)$ 
	which provides us some decay for large times since $k \in \Z^d \setminus \{0\}$. 
	We compute explicitly the integral inside it:
		$$
		 \int_0^t \chi^2_\eps(\sigma) \,  d\sigma = {t \over \eps^2} + {1 \over {2\eps^3}} \left(-e^{-2\eps t}+ 4 e^{-\eps t} -3 \right).
		$$
	Since $e^{-\eps t} \in (0,1]$ for any $t \ge 0$, $-y^2 + 4y -3\in [-3,0]$ for any $y \in [0,1]$ and $|k| \ge 1$, we then write the following bound:
		$$
		\forall \, t \ge 3/(2\eps), \quad \int_0^t \chi^2_\eps(\sigma) \, |k|^2 \, d\sigma \ge {1 \over \eps^2} \left(t - {3 \over {2\eps}}\right) |k|^2 \ge {1 \over \eps^2} \left(t - {3 \over {2\eps}}\right) |k|. 
		$$
	As a consequence, we get
		$$
		\forall \, t \ge 3/\eps, \quad \hbox{exp} \left( - \eps \int_0^t \chi^2_\eps(\sigma) \, |k|^2 \, d\sigma \right) \le e^{-{|k| \over {\eps}} \left(t-{3 \over {2\eps}}\right)} \le e^{-|k| t/(2\eps)}.
		$$
	from which we have:
		\beqn \label{FPeq:exp}
		\forall \, t \ge 3/\eps, \quad e^{\eps t} \, \hbox{exp} \left( - \eps \int_0^t \chi^2_\eps(\sigma) \, |k|^2 \, d\sigma \right) \le e^{-{|k|t /2}} \le e^{-t/2}  
		\eeqn
	where we used that $\eps \in (0,1/2]$.
	Then, 
		$$
		\bal
			I_3 &\le {C \over |k|}  \, \int_{3/\eps}^\infty e^{- t/2} \left| k \, \chi_\eps (t) \,M(k\chi_\eps(t)) - k \, t \, M(kt) \right|\, e^{\zeta t} \, dt \\
			&\le {C \over |k|} \, \int_{3/\eps}^\infty e^{- t/2}\, |k| \, |\chi_\eps(t) - t| \, dt \\
			&\le {{C\eps}} \,   \int_{3/\eps}^\infty t^2 \, e^{- t/2} \, dt
		\eal
		$$
	where we used that the derivative of $y \mapsto y \, M(y)$ (where we still denote the function $M: y \in \R \mapsto (2\pi)^{-3/2} e^{-y^2/2}$) is bounded on $\R$ and the inequality~\eqref{FPeq:chieps}. 
	It implies 
		\beqn \label{FPeq:I3} 
			I_3 \le {{C\eps} \over |k|^3} \,   \int_{0}^\infty t^2 \, e^{- t/2} \, dt \le C \, \eps. 
		\eeqn}

	{ To deal with $I_4$, we use~\eqref{FPeq:exp} and the fact that $\mathds{1}_{t \ge 3/\eps} \le \eps \, t/3$, we thus deduce that 
		\beqn \label{FPeq:I4}
		I_4 \le C \, \int_{0}^\infty t \, M(kt) \, e^{\zeta t} \, \mathds{1}_{t \ge 3/\eps} \, dt \le  {{C\, \eps} \over |k|^3}  \int_0 ^\infty t^2 \, e^{-t^2/2} \, dt  \le C \, \eps. 
		\eeqn}
	
	Gathering~\eqref{FPeq:I11},~\eqref{FPeq:I12},~\eqref{FPeq:I21},~\eqref{Fpeq:I22},~\eqref{FPeq:I3} and~\eqref{FPeq:I4} yields
		\beqn \label{FPeq:Keps-K0}
			| \widetilde{K}_\eps(\tau,k) - \widetilde{K}_0(\tau,k)| \le c_0 \, \eps,
		\eeqn
	for some constant $c_0>0$.
	We then combine~\eqref{FPeq:K0} with~\eqref{FPeq:Keps-K0} to conclude that 
		$$
		\forall \, \tau \in \C, \, \, \hbox{Im} \tau \le 0, \quad  |\widetilde{K_\eps}(\tau,k) - 1| \ge \kappa_0 - c_0 \, \eps.
		$$
	We choose $\eps_0 \in (0,1/12]$ small enough such that $ \kappa_0 - c_0 \, \eps_0>0$, set $\kappa := \kappa_0 - c_0 \, \eps_0$ and it follows that {\bf (H)} is satisfied. 
\qed

\medskip
\subsection{The homogeneous Fokker-Planck equation}
In this part, we recall some features of solutions to the homogeneous Fokker-Planck equation 
		\beqn \label{FPeq:homFP}
		\left\{
			\bal 
			&\partial_t f = \eps \, \left(\Delta_v f + \hbox{div}_v(vf)\right) \\
			&f(0,\cdot) = f_0
			\eal
		\right.
		\eeqn
that will be useful for the study of the evolution of $\widehat{h}(t,0,\xi)$, $t \ge 0$ and $\xi \in \R^d$ in the proof of estimate~\eqref{FPeq:finalh}.
	
We introduce  $L^2_v(m)$ the weighted Lebesgue space associated to the norm 
	\beqn \label{FPeq:L2norm}
	\|f \|_{L^2_v(m)}^2 := \int_{\R^d} f^2(v) \, m^2(v) \, dv, \quad m(v) = \la v \ra^\ell
	\eeqn
where we recall that $\ell$ has been fixed in Subsection~\ref{subsec:not} and satisfies $\ell>d/2$. 
	
The behaviour of solutions to the homogeneous Fokker-Planck equation has been largely studied in several types of functional spaces thanks to Poincar\'{e} inequality or Log-Sobolev inequality for example. 
In our case, we will make use of a result from Gualdani, Mischler and Mouhot~\cite[Theorem~3.1]{GMM} which is stated in the following lemma, 
the main idea is that they have enlarged the space in which the solutions to the equation decay as time goes to infinity to various kind of weighted Lebesgue spaces and in particular to the space $L^2_v(m)$ defined above. 
	
\begin{lem} \label{FPlem:homFP}
	Consider $f_0=f_0(v) \in L^2_v(m)$ such that $\int_{\R^d} f_0(v) \, dv =0$. Then, there exist $C \ge 1$ and $\lambda_*>0$ such that the solution $f=f(t,v)$ to the homogeneous Fokker-Planck equation~\eqref{FPeq:homFP} satisfies
		$$
		\|f(t,\cdot)\|_{L^2_v(m)} \le C \, e^{-\eps \lambda_* t} \, \|f_0\|_{L^2_v(m)}.
		$$
\end{lem}

\medskip
\subsection{Proof of Theorem~\ref{theo1}}
We now handle the proof of Theorem~\ref{theo1} and we split it into two parts, we first look at the behaviour of the density $\rho = \rho(t,x)$ before going back to the solution $h=h(t,x,v)$ itself. 
Let us underline the fact that in the proof of Theorem~\ref{theo0}, we are able to deal with the cases $\eps=0$ and $\eps>0$ in the same time, which is not the case here. 
However, we do not write the proof for $\eps=0$ since it is just a (simpler) adaptation of the case $\eps>0$.
\smallskip

\noindent {\it Proof of the estimate~\eqref{FPeq:finalrho}.}
	As in the proof of the estimate~\eqref{B1eq:finalrho}, we start by noticing that, 
	due to the conservation properties of the operator $- v \cdot \nabla_x h - F(t,x) \cdot \nabla_v M + \eps \, (\Delta_v h + \hbox{div}_v (vh))$, 
	we know that for all times, $\widehat{\rho}(t,0) = (2 \pi)^{d/2} \, \widehat{h_{in}}(0,0)=0$ since $h$ is mean-zero. 
	
	We then concentrate on the study of $\widehat{\rho}(t,k)$ for $k \in \Z^d \setminus \{0\}$. In this part, we consider $\eps \in (0,\eps_0]$ where $\eps_0$ is given in Lemma~\ref{FPlem:Keps} and $t \in \R^+$.
	We first take the Fourier transform of the equation~\eqref{FPeq:VPFP}. It provides us the following equality:
		$$
		\partial_t \widehat{h} - k \cdot \nabla_\xi \widehat{h} + k \cdot \xi \, \widehat{W}(k) \, \widehat{\rho}(t,k) \, \widehat{M}(\xi) + \eps \,  \xi \cdot \nabla_\xi \widehat{h} = - \eps \,  |\xi|^2 \widehat{h} 
		$$
	where we have denoted $\widehat{h} = \widehat{h}(t,k,\xi)$. We rewrite the latter equality as:
		$$
		\partial_t \widehat{h} +(\eps \, \xi -k) \cdot \nabla_\xi \widehat{h} = - \eps \,  |\xi|^2 \widehat{h} - k \cdot \xi \, \widehat{W}(k) \, \widehat{\rho}(t,k) \, \widehat{M}(\xi)
		$$
	in order to consider $(\eps \, \xi -k) \cdot \nabla_\xi$ as our ``new'' transport operator, the term $- \eps \,  |\xi|^2 \widehat{h}$ is seen as an absorption term 
	and $- k \cdot \xi \, \widehat{W}(k) \, \widehat{\rho}(t,k) \, \widehat{M}(\xi)$ as a source one. Then using the method of characteristics, we end up with the following formula:
		\beqn \label{FPeq:Fourier1}
		\bal
			&\widehat{h}(t,k,\xi) = \hbox{exp}\left(- \eps \int_0^t \left| e^{-\eps(t-s)} \xi +\chi_\eps(t-s) \, k \right|^2 \, ds\right) \, \widehat{h_{in}} \left( k, e^{-\eps t} \xi + \chi_\eps(t) \, k \right) \\
			&\quad - \int_0^t  \hbox{exp}\left(- \eps \int_s^t \left| e^{-\eps(t-\sigma)} \xi +
			\chi_\eps(t- \sigma) \,  k \right|^2 \, d\sigma \right) \, k \cdot \left( e^{-\eps(t-s)} \xi + \chi_\eps(t-s) \, k\right) \, \\
			 &\qquad \qquad \qquad \widehat{W}(k) \, \widehat{\rho}(s,k) \, \widehat{M}\left(e^{-\eps(t-s)} \xi +  \chi_\eps(t-s) \, k\right) \, ds.
		\eal
		\eeqn
	Taking $\xi =0$ in the previous formula gives
		$$
		\bal
			&\widehat{\rho}(t,k) =  \hbox{exp}\left(- \eps \int_0^t \left|\chi_\eps(t-s) \, k \right|^2 \, ds\right) \, \widehat{h_{in}} \left( k, \chi_\eps(t) \, k \right) \\
			&\quad - \int_0^t  \hbox{exp}\left(-\eps \int_s^t\left|\chi_\eps(t-\sigma) \, k \right|^2 \, d\sigma \right)  \, k \cdot \left(\chi_\eps(t-s) \, k\right) \widehat{W}(k) \, \widehat{\rho}(s,k) \, \widehat{M}\left( \chi_\eps(t-s) \, k\right) \, ds.
		\eal
		$$
	and thus
		{ 
		$$
		e^{\eps t} \, \widehat{\rho}(t,k) 
		= e^{\eps t} \, \hbox{exp}\left(- \eps \int_0^t \left| \chi_\eps(t-s) \, k \right|^2 \, ds\right) \,\widehat{h_{in}} \left( k, \chi_\eps(t) \, k \right) + \int_0^t K_\eps(t-s) \, e^{\eps s} \, \widehat{\rho}(s,k) \, ds
		$$
	where we recall that $K_\eps$ is defined in~\eqref{FPdef:Keps}.
	Then, using Lemma~\ref{FPlem:Keps} and proceeding as in the proof of estimate~\eqref{B1eq:finalrho}, we obtain for any $k \in \Z^d \setminus \{0\}$ and any $T>0$:
		$$
		\bal
			\sup_{t \in [0,T]} \la kt \ra^n \, e^{\eps t} \, |\widehat{\rho}(t,k)| 
			&\le  \sup_{t \in [0,T]}   \la kt \ra^n \, | \widehat{h_{in}} \left( k, \chi_\eps(t) \, k \right)| \, e^{\eps t} \, \hbox{exp}\left(- \eps \int_0^t \left| \chi_\eps(t-s) \, k \right|^2 \, ds\right)\\
			&\le \sup_{t \in [0,T]}  \left( \frac{\la kt \ra^n}{\la \chi_\eps(t) k \ra^n} \mathds{1}_{0 \le t \le 3/\eps} \, e^3+ \la kt \ra^n e^{-|k| t/2} \mathds{1}_{t \ge 3/\eps}\right)
			\|h_{in}\|_{L^2_x \HH^n_v}
		\eal
		$$
	where we used the estimate~\eqref{FPeq:exp}. 
	But, we recall that $\chi_\eps(t) \ge e^{-3}  t$ for any $t \in [0,3/\eps]$. Consequently
		\beqn \label{FPeq:rho}
			\forall \, k \in \Z^d \setminus \{0\}, \quad \sup_{t \in [0,T]} \la kt \ra^n \, e^{\eps t} \, |\widehat{\rho}(t,k)| 
			\le C(n) \|h_{in}\|_{L^2_x\HH^n_v},
		\eeqn
	this provides us the result~\eqref{FPeq:finalrho}.}
\qed

\smallskip
We are now able to go back to the analysis of the behaviour of the solution $h=h(t,x,v)$ thanks to the previous study on the density $\rho=\rho(t,x)$. 

\smallskip

\noindent {\it Proof of the estimate~\eqref{FPeq:finalh}.} 
	We first deal with the case $k=0$ and consider $\eps \ge 0$. 
	One can notice that the moment of order $0$ in $x$ is not preserved due to the presence of the collision operator as in the Proof of estimate~\eqref{B1eq:finalh}. 
	Indeed, $\la h \ra := \int_{\T^d} h(t,x, \cdot) \, dx$ satisfies the following homogeneous Fokker-Planck equation:
		$$
		\partial_t \la h \ra = \eps \left( \Delta_v \la h \ra + \hbox{div}_v( v \la h \ra)\right).
		$$
	Note that since $h_{in} \in L^2_x\HH^n_v$, we have $\la h_{in}\ra \in L^2_v(m)$ where $L^2_v(m)$ is defined through its norm~\eqref{FPeq:L2norm}. Indeed, using Jensen inequality,
		$$
		\bal
			&\int_{\R^d} \la h_{in} \ra^2 \, m^2(v) \, dv = \int_{\R^d} \left( \int_{\T^d} h_{in}(x,v) \, dx \right)^2 \, m^2(v) \, dv \\
			&\qquad \qquad \le C \, \int_{\R^d \times \T^d} h^2_{in}(x,v) \, \la v \ra^{2 \ell} \,dx \, dv \le C \, \|h_{in}\|^2_{L^2_x \HH^n_v}.
		\eal
		$$
	Using Lemma~\ref{FPlem:homFP}, we have
		$$
		\forall \, t \ge 0, \quad \| \la h \ra(t,\cdot) \|_{L^2_v(m)} \le C \, e^{- \eps \lambda_* t} \| \la h_{in} \ra \|_{L^2_v(m)}
		$$
	where we used that $h_{in}$ is mean-zero. Consequently, using Cauchy-Schwarz inequality combined with the fact that $m^{-2} = \la \cdot \ra^{-2 \ell} \in L^1_v(\R^d)$, we have 
		$$
		\bal
			&\forall \, \xi \in \R^d,\, \forall \, t \ge 0, \\ 
			&\qquad \qquad |\widehat{h}(t,0,\xi)| = |\widehat{\la h \ra}(t,\xi)| = \left| \int_{\R^d} \la h \ra(t,v) \, m(v) \, m^{-1}(v) \, e^{-i v \cdot \xi} \, dv \right| \\
			&\qquad \qquad \qquad \le \| \la h \ra (t, \cdot) \|_{L^2_v(m)} \, \int_{\R^d} m^{-2}(v) \, dv \\
			&\qquad \qquad \qquad\le C \, e^{- \eps \lambda_* t} \| \la h_{in} \ra \|_{L^2_v(m)} \le C \, e^{- \eps \lambda_* t} \| h_{in} \|_{L^2_x \HH^n_v}
		\eal
		$$
	from which we readily conclude.

\smallskip
	{ We now deal with $k \in \Z^d \setminus \{0\}$ and consider $\eps \in (0,\eps_0]$ where $\eps_0$ is given in Lemma~\ref{FPlem:Keps}.
	Going back to the equality~\eqref{FPeq:Fourier1}, we have:
		\beqn \label{FPeq:Fourierh2}
		\bal
			&e^{\eps t} \,\widehat{h}(t,k,\xi) = e^{\eps t} \, \hbox{exp}\left(- \eps \int_0^t \left| e^{-\eps s} \xi+\chi_\eps(s) k \right|^2 \, ds\right) \, \widehat{h_{in}} \left( k, e^{-\eps t}\xi + \chi_\eps(t)k \right) \\
			&\quad - \int_0^t  e^{\eps s} \,\hbox{exp}\left(- \eps \int_0^s  \left| e^{-\eps \sigma} \xi +\chi_\eps(\sigma) k \right|^2 \, d\sigma \right) \, k \cdot  \left( e^{-\eps s} \xi +\chi_\eps(s) k \right) \\
			&\qquad \qquad \qquad \widehat{W}(k) \, e^{\eps (t-s)} \,\widehat{\rho}(t-s,k) \, \widehat{M}\left(e^{-\eps s} \xi +\chi_\eps(s) k\right) \, ds =: T_1+T_2.
		\eal
		\eeqn
	We start by studying the terms of type $\hbox{exp}\left(- \eps \int_0^t \left| e^{-\eps s} \xi+\chi_\eps(s) k \right|^2 \, ds\right)$. 
	As previously, we compute explicitly the integral inside it:
		$$
		\bal
			\eps \, \int_0^t \left| e^{-\eps s} \xi+\chi_\eps(s) k \right|^2 \, ds 
			&= {{1-e^{-2\eps t}} \over 2} |\xi|^2 + {{(1-e^{-\eps t})^2} \over \eps} \, \xi \cdot k \\
			&\qquad + \left({t \over \eps} + {1 \over {2\eps^2}} \left[-e^{-2\eps t} + 4 e^{-\eps t} -3 \right] \right) |k|^2.
		\eal
		$$
	Then using the following bound:
		$$
		{{(1-e^{-\eps t})^2} \over \eps} \, |\xi \cdot k| \le {|\xi \cdot k| \over \eps} \le {|k|^2 \over {\eps^2}} + {|\xi|^2 \over 4},
		$$
	we deduce:
		$$
		\bal
			\eps \int_0^t \left| e^{-\eps s} \xi+\chi_\eps(s) k \right|^2 \, ds 
			&\ge{{1-2e^{-2\eps t}} \over 4} |\xi|^2 	+ \left({t \over \eps} + {1 \over {2\eps^2}} \left[-e^{-2\eps t} + 4 e^{-\eps t} -5 \right] \right) |k|^2 \\
			&\ge {1 \over \eps} \left(t-{5\over \eps}\right)|k|^2
		\eal
		$$
	where the last inequality holds for $t \ge \ln(2)/(2\eps)$ and in particular for $t \ge 6/\eps$ and where we used that $-y^2+4y-5 \ge -5$ for any $y \in [0,1]$. 
	From this, we have:
		$$
			\forall \, t \ge 6/\eps,\quad  \hbox{exp}\left(- \eps \int_0^t \left| e^{-\eps s} \xi+\chi_\eps(s) k \right|^2 \, ds\right) \le e^{-t/(6\eps)} 
		$$
	and thus 
		\beqn \label{FPeq:exp3}
			\forall \, t \ge 6/\eps,\quad e^{\eps t} \, \hbox{exp}\left(- \eps \int_0^t \left| e^{-\eps s} \xi+\chi_\eps(s) k \right|^2 \, ds\right) \le e^{\eps t} \, e^{-t/(6\eps)} \le e^{-t/6} 
		\eeqn
	where we used that $\eps \in (0, 1/3]$. 
	Let us now estimate 	each term of~\eqref{FPeq:Fourierh2}. If $t \ge 6/\eps$, the first term is easily treated using~\eqref{FPeq:exp3}:
		\beqn \label{FPeq:T11}
			|T_1| \le C \, e^{-t/6} \le C \, \langle t \rangle^{2-n}.
		\eeqn
	If $0 \le t \le 6/\eps$, we have 
		\beqn \label{FPeq:T12}
			\bal
				|T_1| &\le C \, \la e^{-\eps t}\xi + \chi_\eps(t)k \ra^{-n} \\
				&\le C \, \la \xi \ra^{2-n} \la \chi_\eps(t) k \ra^{2-n} \\
				&\le C \, \la \xi \ra^{2-n} \la t \ra^{2-n}
			\eal 
		\eeqn
	where we used that $\chi_\eps(t) \ge e^{-6} t$ for $0 \le t \le 6/\eps$.
	Concerning the second one, we decompose it into two parts:
		$$
		\bal
			&T_2 =   \int_0^t  \mathds{1}_{s\le 6/\eps} \,  e^{\eps s} \, 
			\hbox{exp}\left(- \eps \int_0^s  \left| e^{-\eps \sigma} \xi +\chi_\eps(\sigma) k \right|^2 \, d\sigma \right) \, k \cdot  \left( e^{-\eps s} \xi +\chi_\eps(s) k \right) \\
			&\qquad \qquad \qquad \widehat{W}(k) \,e^{\eps (t-s)} \,  \widehat{\rho}(t-s,k) \, {M}\left(e^{-\eps s} \xi +\chi_\eps(s) k\right) \, ds \\
			&\qquad +  \int_0^t  \mathds{1}_{s\ge 6/\eps} \,  e^{\eps s} \, 
			\hbox{exp}\left(- \eps \int_0^s  \left| e^{-\eps \sigma} \xi +\chi_\eps(\sigma) k \right|^2 \, d\sigma \right) \, k \cdot  \left( e^{-\eps s} \xi +\chi_\eps(s) k \right) \\
			&\qquad \qquad \qquad \widehat{W}(k) \, e^{\eps (t-s)} \, \widehat{\rho}(t-s,k) \, {M}\left(e^{-\eps s} \xi +\chi_\eps(s) k\right) \, ds \\
			&=:T_{21}+T_{22}. 
		\eal
		$$
	To handle the first term $T_{21}$, we use the decay property of $\widehat{M} = M$ which in particularly implies:
		$$
		\left| e^{-\eps s} \xi -\chi_\eps(s) k \right| \widehat{M}\left(e^{-\eps s} \xi -\chi_\eps(s) k\right) \le C \, \la e^{-\eps s} \xi -\chi_\eps(s) k \ra^{2-n}.
		$$
	and the estimate~\eqref{FPeq:finalrho} so that for any $0 \le s \le 6/\eps$:
		$$
		\bal
			&e^{\eps (t-s)} \,|\widehat{\rho}(t-s,k)| \left| e^{-\eps s} \xi -\chi_\eps(s) k \right| \widehat{M}\left(e^{-\eps s} \xi -\chi_\eps(s) k\right) \\
			&\qquad \le C \, \langle t-s \rangle^{2-n} \, \langle t-s \rangle^{-2}  \, \la e^{-\eps s} \xi -\chi_\eps(s) k \ra^{2-n} \\
			&\qquad \le C \, \langle t \rangle^{2-n} \, \langle s \rangle^{n-2} \, \langle t-s \rangle^{-2} \, \la \xi \ra^{n-2} \la \chi_\eps(s) k \ra^{2-n} \\
			&\qquad \le C \, \langle t \rangle^{2-n} \, \la \xi \ra^{n-2}\, \langle t-s \rangle^{-2}
		\eal
		$$
	where we used that $\chi_\eps(t) \ge e^{-6} s$ for $0 \le s \le 6/\eps$ to get the last inequality. As a consequence, we obtain:
		\beqn \label{FPeq:T21}
		\bal
			T_{21} &\le C \, \langle t \rangle^{2-n} \, \la \xi \ra^{n-2} \int_0^{6/\eps} \langle t-s \rangle^{-2} \, ds \\
			&\le C \, \langle t \rangle^{2-n} \, \la \xi \ra^{n-2} \int_0^\infty \la s \ra^{-2} \, ds \\
			&\le C \, \langle t \rangle^{2-n} \, \la \xi \ra^{n-2}.
		\eal
		\eeqn
	The second term is treated thanks to~\eqref{FPeq:exp3},~\eqref{FPeq:finalrho} and the fact that $y \mapsto |y| M(y)$ is bounded on $\R^d$:
		\beqn \label{FPeq:T22}
		\bal
			T_{22} &\le C \int_0^t e^{-s/6} \, \la t-s \ra^{-n} \, ds \\
			&\le  C \, \la t \ra^{-n} \int_0^\infty e^{-s/6} \, \la s \ra^{n} \, ds \\
			&\le C  \, \la t \ra^{-n}.
		\eal
		\eeqn
	Gathering~\eqref{FPeq:T11},~\eqref{FPeq:T12},~\eqref{FPeq:T21},~\eqref{FPeq:T22}, we conclude that~\eqref{FPeq:finalh} holds. }
\qed
\smallskip

%




\bigskip
\bigskip
\bibliographystyle{acm}

\end{document}